\documentclass[11pt]{article}
\renewcommand{\baselinestretch}{1.2}
\newcommand{\dated}{\mbox{} \hfill {\small [{\tt \today}]}} 
\usepackage{amsthm,enumerate}
\theoremstyle{plain}
\newtheorem{theorem}{Theorem}[section]
\newtheorem{lemma}[theorem]{Lemma}

\newtheorem{proposition}[theorem]{Proposition}
\theoremstyle{definition}
\newtheorem{definition}[theorem]{Definition}
\theoremstyle{remark}
\newtheorem*{remark}{Remark}
\newtheorem*{example}{Example}
\newtheorem*{rems}{Remarks}
\newtheorem*{exs}{Examples}
\newenvironment{remarks}{\begin{rems}\begin{enumerate}}{\end{enumerate}\end{rems}}

\newenvironment{items}{\begin{enumerate}[\rm (i)]}{\end{enumerate}}
\newenvironment{alphitems}{\begin{enumerate}[\rm (a)]}{\end{enumerate}}

 \usepackage{amsmath,amssymb,amsfonts,diagrams}
\newenvironment{keywords}{\noindent\small {\it Keywords\/}:}{\vskip 4pt}
\newenvironment{classification}{\noindent\small 2000 {\it Mathematics Subject
Classification\/}:}{\vskip 12pt}

\newcommand{\comps}{{\mathbb C}}

\newcommand{\tensor}{\otimes}
\newcommand{\ttensor}{\tilde{\otimes}}
\newcommand{\Tensor}{\hat{\otimes}}

\newcommand{\wTensor}{\check{\otimes}}

\newcommand{\cstar}{{C^\ast}}

\newcommand{\id}{{\mathrm{id}}}

\newcommand{\clspan}{{\overline{\operatorname{span}}}}
\newcommand{\lspan}{{\operatorname{span}}}

\newcommand{\A}{{\mathfrak A}}

\newcommand{\Hilbert}{{\mathfrak H}}

\newcommand{\M}{{\mathfrak M}}

\newcommand{\CB}{{\cal CB}}
\newcommand{\op}{{\mathrm{op}}}

\newcommand{\VN}{\operatorname{VN}}

\newcommand{\varcl}[1]{\overline{#1}}

\newcommand{\LUC}{\mathit{LUC}}
\newcommand{\RUC}{\mathit{RUC}}
\newcommand{\UC}{\mathit{UC}}
\newcommand{\WAP}{\mathit{WAP}}
\newcommand{\G}{\mathbb{G}}
\title{Uniform continuity over locally compact quantum groups}
\author{\textit{Volker Runde}}
\date{}
\begin{document}
\maketitle
\begin{abstract}
We define, for a locally compact quantum group $\G$ in the sense of Kustermans--Vaes, the space of $\LUC(\G)$ of left uniformly continuous elements in $L^\infty(\G)$. This definition covers both the usual left uniformly continuous functions on a locally compact group and Granirer's uniformly continuous functionals on the Fourier algebra. We show that $\LUC(\G)$ is an operator system containing the $\cstar$-algebra ${\cal C}_0(\G)$ and contained in its multiplier algebra ${\cal M}({\cal C}_0(\G))$. We use this to partially answer an open problem by B\'edos--Tuset: if $\G$ is co-amenable, then the existence of a left invariant mean on ${\cal M}({\cal C}_0(\G))$ is sufficient for $\G$ to be amenable. Furthermore, we study the space $\WAP(\G)$ of weakly almost periodic elements of $L^\infty(\G)$: it is a closed operator system in $L^\infty(\G)$ containing ${\cal C}_0(\G)$ and---for co-amenable $\G$---contained in $\LUC(\G)$. Finally, we show that---under certain conditions, which are always satisfied if $\G$ is a group---the operator system $\LUC(\G)$ is a $\cstar$-algebra.
\end{abstract}
\section*{Introduction}
For any locally compact group $G$, and a function $f \!: G \to \comps$, we denote by $L_x f \!: G \to \comps$ the \emph{left translate} of $f$ by $x \in G$, i.e., $(L_x f)(y) = f(xy)$ for $y \in G$. If $f \in {\cal C}(G)$, i.e., is bounded and continuous, we call it \emph{left uniformly continuous} if the map
\[
  G \to {\cal C}(G), \quad x \mapsto L_x f 
\]
is continuous with respect to the given topology on $G$ and the norm topology on ${\cal C}(G)$. The collection of all left uniformly continuous functions on $G$---denoted by $\LUC(G)$---is obviously a unital $\cstar$-subalgebra of ${\cal C}(G)$.
Somewhat less obvious is the fact that $\LUC(G)$ consists precisely of the functions $\phi \cdot f$ with $\phi \in L^\infty(G)$ and $f \in L^1(G)$, where $\cdot$ denotes the canonical module action of the Banach algebra $L^1(G)$ on its dual space $L^\infty(G)$ (\cite[(32.45)(a) and (b)]{HR}; in \cite{HR}, a left uniformly continuous function in our sense is called right uniformly continuous.) Similarly, one defines \emph{right uniformly continuous} functions on $G$---denoted by $\RUC(G)$---, and calls the functions in $\UC(G) := \LUC(G) \cap \RUC(G)$ \emph{uniformly continuous}. Left and right uniform continuity are important concepts in the study of locally compact groups. For instance, $\LUC(G) = \RUC(G)$ holds if and only if $G$ has small invariant neighborhoods (\cite[(4.14)(g)]{HR}), and this is the case if and only if $L^1(G)$ has a bounded approximate identity in its center (\cite{Mos}). The spaces $\LUC(G)$, $\RUC(G)$, and $\UC(G)$, also play an important r\^ole in the theory of amenable groups as natural domains for invariant means (see \cite{Pie}).
\par 
In \cite{Eym}, P.\ Eymard defined the Fourier algebra $A(G)$ for an arbitrary locally compact group. Like $L^1(G)$, it is the predual of a von Neumann algebra, namely of $\VN(G)$, the group von Neumann algebra of $G$, which is generated by the left regular representation of $G$ on $L^2(G)$. In \cite{Gra1}, E.\ E.\ Granirer defined the space $\UC(\hat{G})$ uniformly continuous functionals on $A(G)$ as the closed linear span of $\{ x \cdot f : x \in \VN(G),\, f \in A(G) \}$, where $\cdot$ stands for the canonical action of $A(G)$ on its dual $\VN(G)$. Even though it is not obvious from this definition, $\UC(\hat{G})$ is indeed a $\cstar$-subalgebra of $\VN(G)$: it is the norm closure of the operators is in $\VN(G)$ with compact support in the sense of \cite{Eym}. It contains $C^\ast_r(G)$, the reduced group $\cstar$-algebra of $G$ (\cite[Proposition 2]{Gra2}), and is contained in its multiplier algebra ${\cal M}(C^\ast_r(G))$ (\cite[Proposition 1]{Gra2}).
\par 
In the past decades, various attempts have been made to develop a rigorous framework for a duality theory for general locally compact groups that extends the Pontryagin duality for locally compact abelian groups. One such framework is the theory of Kac algebras, as expounded in the monograph \cite{ES}. In the Kac algebra framework, the dual of $L^1(G)$ is $A(G)$ in a precise and well defined manner, and in view of the parallelism between uniform continuity in $L^\infty(G)$ and in $\VN(G)$, one is tempted to develop a unified notion of uniform continuity for both in a general Kac algebraic setting. One drawback of such an endeavor, however, is the lack of examples for Kac algebras: \cite{ES} only gives two kinds of examples, namely those correspoding to locally compact groups---via $L^\infty(G)$---and those---via $\VN(G)$---dual to them.
\par 
Fairly recently, J.\ Kustermans and S.\ Vaes proposed a surprisingly simple set of axioms for so-called locally compact quantum groups (\cite{KV1} and \cite{KV2}). For a detailed exposition on the evolution of these axioms---with many references to the original literature---, we refer to the introduction of \cite{KV1} and to \cite{Vai}. The Kustermans--Vaes axioms cover the Kac algebras (and therefore all locally compact groups), allow for the development of a Pontryagin type duality theory, but also seem to be satisfied by all known examples of $\cstar$-algebraic quantum groups, such as Woronowicz's $\mathrm{SU}_q(2)$ (\cite{Wor}). For a simultaneous treatment of uniform continuity both in $L^\infty(G)$ and in $\VN(G)$, the framework of locally compact quantum groups thus seems to be best suited.
\par 
In this paper, we define, for a locally compact quantum group $\G$ (the notation will be explained in Section \ref{prelim} below), the space $\LUC(\G)$ of left uniformly continuous elements in $L^\infty(\G)$ (and, analogously, $\RUC(\G)$ and $\UC(\G)$). These definitions, of which $\RUC(\G)$ already appears in \cite{HNR}, simultaneously cover both $\LUC(G)$, etc., as well as $\UC(\hat{G})$ for a locally compact group $G$. We show:
\begin{itemize}
\item $\LUC(\G)$ is an operator system containing the $\cstar$-algebra ${\cal C}_0(\G)$ and contained in its multiplier algebra ${\cal M}({\cal C}_0(\G))$.
\item For co-amenable $\G$, the amenability of $\G$ is already ensured by the existence of a left invariant mean on $\LUC(\G)$ or ${\cal M}({\cal C}_0(\G))$; this partially answers an open problem brought up by E.\ B\'edos and L.\ Tuset (\cite[p.\ 876]{BT}).
\item The space $\WAP(\G)$ of weakly almost periodic elements of $L^\infty(\G)$ contains ${\cal C}_0(\G)$ and, for co-amenable $\G$, is contained in $\UC(\G)$.
\item Under certain conditions, which are always satisfied if $\G$ is a group, the operator system $\LUC(\G)$ is, in fact, a $\cstar$-algebra.
\end{itemize}
\section{Locally compact quantum groups---an overview} \label{prelim}
In this preliminary section, we give a brief overview of locally compact quantum groups in the sense of Kustermans and Vaes (\cite{KV1} and \cite{KV2}). We focus on the von Neumann algebraic approach, as expounded in \cite{KV2} or \cite{vDa}, where the latter reference presents von Neumann algebraic quantum groups independent of the $\cstar$-algebraic approach in \cite{KV1}. Nevertheless, we shall require some facts about (reduced) $\cstar$-algebraic quantum groups as well. For details, we refer to \cite{KV1}, \cite{KV2}, and \cite{vDa}. We shall also---not so much in this section, but later on---require results from the theory of operator spaces. For background on this theory, we use \cite{ER} as a reference and adopt that book's notation; in particular, $\Tensor$ and $\wTensor$ stand for the injective and projective tensor product, respectively, of operator spaces and not of Banach spaces. (Restricted to $\cstar$-algebras, $\wTensor$ is just the spatial tensor product.)
\par 
As a (von Neumann algebraic) locally compact quantum group is a Hopf--von Neumann algebra with additional structure, we begin with recalling the definition of a Hopf--von Neumann algebra ($\bar{\tensor}$ denotes the $W^\ast$-tensor product):
\begin{definition}
A \emph{Hopf--von Neumann algebra} is a pair $(\M,\Gamma)$, where $\M$ is a von Neumann algebra and $\Gamma \!: \M \to \M \bar{\tensor} \M$ is a \emph{co-multiplication}, i.e., a normal, unital $^\ast$-homomorphism satisfying $(\id \tensor \Gamma) \circ \Gamma = (\Gamma \tensor \id) \circ \Gamma$.
\end{definition}
\begin{example}
For a locally compact group $G$, define $\Gamma_G \!: L^\infty(G) \to L^\infty(G \times G)$ and $\hat{\Gamma}_G \!: \VN(G) \to \VN(G \times G)$ via 
\begin{align*}
  (\Gamma_G \phi)(x,y) & := \phi(xy) \qquad (\phi \in L^\infty(G), \, x,y \in G) \\
\intertext{and}
  \hat{\Gamma}_G(\lambda(x)) & := \lambda(x) \tensor \lambda(x) \qquad (x \in G),
\end{align*}
where $\lambda$ is the left regular representation of $G$ on $L^2(G)$. Then $\Gamma_G$ and $\hat{\Gamma}_G$ are co-multiplications, turning $(L^\infty(G),\Gamma_G)$ and $(\VN(G),\hat{\Gamma}_G)$ into Hopf--von Neumann algebras. 
\end{example}
\begin{remarks}
\item One can define a product $\ast$ on $\M_\ast$, the unique predual of $\M$, turning it into a Banach algebra:
\begin{equation} \label{prod}
  \langle x, f \ast g \rangle := \langle \Gamma x, f \tensor g \rangle \qquad (f,g \in \M_\ast, \,   x \in \M).
\end{equation}
For $(L^\infty(G), \Gamma_G)$, where $G$ is a locally compact group, (\ref{prod}) yields the usual convolution product on $L^1(G)$ whereas for $(\VN(G),\hat{\Gamma}_G)$, it gives us pointwise multiplication in $A(G)$.
\item By \cite[Theorem 7.2.4]{ER}, $(\M \bar{\tensor} \M)_\ast \cong \M_\ast \Tensor \M_\ast$ holds. Since $\Gamma \!: \M\to \M \bar{\tensor} \M$ is a normal, i.e., weak$^\ast$ continuous, injective $^\ast$-homomorphism and thus a complete isometry, it is the adjoint of a map $\Gamma_\ast \!: \M_\ast \Tensor \M_\ast \to \M_\ast$, which is a complete quotient map by \cite[Corollary 4.1.9]{ER}; in particular, 
\begin{equation} \label{squaredense}
   \M_\ast = \clspan \{ f \ast g : f , g \in \M_\ast \}
\end{equation}
holds, where $\clspan$ stands for the closed, linear span. 
\end{remarks}
\par 
To define the additional structure that turns a Hopf--von Neumann algebra into a locally compact quantum group, we recall some basic facts about weights (see \cite{Tak}, for instance). 
\par
Let $\M$ be a von Neumann algebra, and let $\M^+$ denote its positive elements. A \emph{weight} on $\M$ is an additive map $\phi \!: \M^+ \to [0,\infty]$ such that $\phi(tx) = t \phi(x)$ for $t \in [0,\infty)$ and $x \in \M^+$. We let
\begin{align*}
  {\cal M}_\phi^+ & := \{ x \in \M^+ : \phi(x) < \infty \}, \\
  {\cal M}_\phi & := \lspan \, {\cal M}_\phi^+, \\
\intertext{and}
  {\cal N}_\phi & := \{ x \in \M : x^\ast x \in {\cal M}_\phi \}.
\end{align*}
Then $\phi$ extends to a linear map on ${\cal M}_\phi$, and ${\cal N}_\phi$ is a left ideal of $\M$. Using the GNS-construction, we obtain a representation $\pi_\phi$ of $\M$ on some Hilbert space $\Hilbert_\phi$; we denote the canonical map from ${\cal N}_\phi$ into $\Hilbert_\phi$ by $\Lambda_\phi$. Moreover, we call $\phi$ \emph{semi-finite} if ${\cal M}_\phi$ is weak$^\ast$ dense in $\M$, \emph{faithful} if $\phi(x) = 0$ for $x \in \M^+$ implies that $x = 0$, and \emph{normal} if $\sup_\alpha \phi(x_\alpha) = \phi\left( \sup_\alpha x_\alpha \right)$ for each bounded, increasing net $( x_\alpha )_\alpha$ in $\M^+$. 
\begin{definition} \label{lcqg}
A (von Neumann algebraic) \emph{locally compact quantum group} is a Hopf--von Neumann algebra $(\M,\Gamma)$ such that:
\begin{alphitems}
\item there is a normal, semifinite, faithful weight $\phi$ on $\M$---a \emph{left Haar weight}---which is left invariant, i.e., satisfies
\[
  \phi((f \tensor \id)(\Gamma x)) = f(1) \phi(x) \qquad   (f \in \M_\ast, \, x \in {\cal M}_\phi);
\]
\item there is a normal, semifinite, faithful weight $\psi$ on $\M$---a \emph{right Haar weight}---which is right invariant, i.e., satisfies
\[
  \phi((\id \tensor f)(\Gamma x)) = f(1) \psi(x) \qquad (f \in \M_\ast, \, x \in {\cal M}_\psi).
\]
\end{alphitems}
\end{definition}
\begin{example}
Let $G$ be a locally compact group. Then the Hopf--von Neumann algebra $(L^\infty(G),\Gamma_G)$ is a locally compact quantum group: $\phi$ and $\psi$ can be chosen as left and right Haar measure, respectively. For $(\VN(G),\hat{\Gamma}_G)$, the Plancherel weight  (\cite[Definition VII.3.2]{Tak}) is both a left and a right Haar weight.
\end{example}
\begin{remark}
Even though only the existence of a left and a right Haar weight, respectively, is presumed, both weights are actually unique up to a positive scalar multiple. 
\end{remark}
\par
An extremely important object associated with every locally compact quantum group $(\M,\Gamma)$ is its \emph{multiplicative unitary}: it is the unique operator $W \in {\cal B}(\Hilbert_\phi \ttensor_2 \Hilbert_\phi)$, where $\ttensor_2$ stands for the Hilbert space tensor product, satisfying
\begin{equation} \label{fundi*}
  W^\ast(\Lambda_\phi(x) \tensor \Lambda_\phi(y))
  = (\Lambda_\phi \tensor \Lambda_\phi)((\Gamma y)(x \tensor 1))
  \qquad (x,y \in {\cal N}_\phi).
\end{equation}
Using the left invariance of $\phi$, it is easy to see that $W^\ast$ is an isometry whereas it is considerably more difficult to show that $W$ is indeed a unitary operator (\cite[Theorem 3.16]{KV1}). The unitary $W$ lies in $\M \bar{\tensor} {\cal B}(\Hilbert_\phi)$ and implements the co-multiplication via 
\[
  \Gamma x = W^\ast (1 \tensor x) W \qquad (x \in \M)
\]
(see the discussion following \cite[Theorem 1.2]{KV2}). 
\begin{example}
If $G$ is a locally compact group, then the multiplicative unitary of $(L^\infty(G),\Gamma_G)$. Then its multiplicative unitary on $L^2(G \times G)$ is given by
\[
  (W\boldsymbol{\xi})(x,y) := \boldsymbol{\xi}(x,x^{-1}y) \qquad (\boldsymbol{\xi} \in L^2(G \times G), \, x,y \in G).
\]
\end{example}
\par
To emphasize the parallels between locally compact quantum groups and groups, we shall use the following notation (which was suggested by Z.-J.\ Ruan and is used in \cite{DR}, \cite{HNR}, \cite{JNR}, and \cite{Run3}): We refer to a locally compact quantum group $(\M,\Gamma)$ by the symbol $\G$ and write: $L^\infty(\G)$ for $\M$, $L^1(\G)$ for $\M_\ast$, and $L^2(\G,\phi)$ for $\Hilbert_\phi$. If $L^\infty(\G) = L^\infty(G)$ for a locally compact group $G$ and $\Gamma = \Gamma_G$, we say that $\G$ actually \emph{is} the group $G$; a locally compact quantum group $\G$ is of this form precisely if $L^\infty(\G)$ is abelian (this follows from \cite[Th\'eor\`eme 2.2]{BS}).
\par 
Of course, given a locally compact quantum group $\G$, there is no intrinsic reason to give preference to the left Haar weight $\phi$ over the right Haar weight $\psi$. If we perform the GNS-construction with respect to the right Haar weight, we obtain a Hilbert space $L^2(\G,\psi)$, and a unitary $V \in {\cal B}(L^2(\G,\psi) \ttensor_2 L^2(\G,\psi))$ defined in a way similar to (\ref{fundi*}). Even though $\phi$ and $\psi$ do not appear to be related in Definition \ref{lcqg}, they are via the \emph{unitary antipode} $R$ of $\G$: we can always suppose that $\psi = \phi \circ R$. This allows to identify $L^2(\G,\phi)$ and $L^2(\G,\psi)$---which shall henceforth be denoted by simply $L^2(\G)$---and to consider both $W$ and $V$ as operators on $L^2(\G) \ttensor_2 L^2(\G)$. (At times, using $V$ instead of $W$ can be useful; see \cite{JNR}). 
\par
Next, we sketch the duality for locally compact quantum groups. It extends Pontryagin duality for locally compact abelian groups.
\par
For a locally compact quantum group $\G$, set
\[
  L^\infty(\hat{\G}) : = 
  \varcl{\{ (f \tensor \id)(W) : f \in L^1(\G) \}}^{\,\text{$\sigma$-strong$^\ast$}},
\]
where $f \tensor \id$ for $f \in L^1(\G)$ is the usual slice map. It can be shown that $L^\infty(\hat{\G})$ is a von Neumann algebra. Let $\sigma$ denote the flip map on $L^2(\G) \ttensor_2 L^2(\G)$, i.e., $\sigma(\xi \tensor \eta) = \eta \tensor \xi$ for $\xi, \eta \in L^2(\G)$, and set $\hat{W} := \sigma W^\ast \sigma$. Then 
\[
  \hat{\Gamma} \!:  L^\infty(\hat{\G}) \to 
  L^\infty(\hat{\G}) \bar{\tensor}  L^\infty(\hat{\G}), \quad
  x \mapsto \hat{W}^\ast(1 \tensor x) \hat{W}.
\]
is a co-multiplication. One can also define a left Haar weight $\hat{\phi}$ and a right Haar weight $\hat{\psi}$ for $(L^\infty(\hat{\G}),\hat{\Gamma})$ turning it into a locally compact quantum group again, the \emph{dual quantum group} of $\G$, which we denote by $\hat{\G}$, and whose multiplicative unitary is $\hat{W}$ as defined above. Finally, a Pontryagin duality theorem holds, i.e., $\Hat{\Hat{\G}} = \G$. 
\begin{example}
Let $G$ be a locally compact group. Then $L^\infty(\hat{G})$ is the $\sigma$-strong$^\ast$ closure of $\lambda(L^1(G))$, i.e., $L^\infty(\hat{G}) = \VN(G)$, the group von Neumann algebra of $G$. Further,the co-multiplication $\hat{\Gamma}_G \!: L^\infty(\hat{G}) \to L^\infty(\hat{G}) \bar{\tensor} L^\infty(\hat{G})$ is just $\hat{\Gamma}_G \!: \VN(G) \to \VN(G \times G)$ introduced earlier, so that $L^1(\hat{G}) = A(G)$. 
\end{example}
\par 
We conclude this section with a look at the relation between von Neumann algebraic and (reduced) $\cstar$-algebraic quantum groups.
\par 
Given a locally compact quantum group $\G$ with multiplicative unitary $W$, we set
\[
  {\cal C}_0(\G) := \varcl{\{ (\id \tensor \omega)(W) : \omega \in {\cal B}(L^2(\G))_\ast \}}^{\| \cdot \|}
\]
Restricting $\Gamma$ to ${\cal C}_0(\G)$ then yields a reduced $\cstar$-algebraic quantum group in the sense of \cite[Definition 4.1]{KV1} (\cite[Proposition 1.6]{KV2}). For a $\cstar$-algebra $\A$, we denote its multiplier algebra by ${\cal M}(\A)$ (and its left, right, and quasi-multipliers by ${\cal LM}(\A)$, ${\cal RM}(\A)$, and ${\cal QM}(\A)$, respectively; see, for instance, \cite{Ped} for the definitions). The restriction of $\Gamma$ to ${\cal C}_0(\G)$ is a non-degenerate $^\ast$-homomorphism into ${\cal M}({\cal C}_0(\G) \wTensor {\cal C}_0(\G))$; by strict continuity, it extends to a $^\ast$-homomorphism from ${\cal M}({\cal C}_0(\G))$ to ${\cal M}({\cal C}_0(\G) \wTensor {\cal C}_0(\G))$. (We would have obtained the same map by restricting $\Gamma$ to ${\cal M}({\cal C}_0(\G))$ right away.) A fact we shall need later on is that the multiplicative unitary $W$ of $\G$ does not just lie in $L^\infty(\G) \bar{\tensor} {\cal B}(L^2(\G))$, but in fact already---with the appropriate identifications in place---in ${\cal M}({\cal C}_0(\G) \wTensor {\cal K}(L^2(\G)))$ (\cite[pp.\ 912--913]{KV1}). The dual space of ${\cal C}_0(\G)$---suggestively denoted by $M(\G)$---becomes also a Banach algebra with a product defined in a way similar to (\ref{prod}); it canonically contains $L^1(\G)$ as a closed ideal (\cite[pp.\ 193--194]{KV1}). We shall denote the product in $M(\G)$ by $\ast$ as well.
\begin{example}
Let $G$ be a locally compact group. Then ${\cal C}_0(G)$ and $M(G)$ in the sense just outlined are the usual objects denoted by those symbols whereas ${\cal C}_0(\hat{G})$ is $C^\ast_r(G)$ and $M(\hat{G})$ is the reduced Fourier--Stieltjes algebra of \cite{Eym}.
\end{example}
\section{$\LUC(\G)$, $\RUC(\G)$, and $\UC(\G)$---definition and basic properties}
Recall (see \cite[Examples 2.6.2(v)]{Dal}), for instance, that any Banach algebra $\A$ and any Banach $\A$-module $E$ (left, right, or bi-) there is a canonical way of turning the dual space $E^\ast$ into a Banach $\A$-module (right, left, or bi-); in particular, this applies to $E = \A$. Given a locally compact quantum group $\G$, the module actions of $L^1(\G)$ on $L^\infty(\G)$ are given by 
\[
  f \cdot x = (\id \tensor f)(\Gamma x) \quad\text{and}\quad x \cdot f = (f \tensor \id)(\Gamma x) 
  \qquad (f \in L^1(\G), \, x \in L^\infty(\G)).
\]
Throughout, we adopt the convention to denote the actions of $L^1(\G)$ on $L^\infty(\G)$, etc., as well the corresponding dual actions by $\cdot$ whereas we express the module actions of $L^\infty(\G)$ on $L^1(\G)$, etc., by mere juxtaposition.
\par 
Before we give our definitions of uniform continuity over a locally compact quantum group, we state and prove the following (somewhat folkloristic) lemma for later reference:
\begin{lemma} \label{faclem}
Let $\A$ be a $\cstar$-algebra, and let $E$ be a closed submodule of the Banach $\A$-bimodule $\A^\ast$. Then we have $E = \{ afb : a,b \in \A, \, f \in E \}$.
\end{lemma}
\begin{proof}
Let $f \in E$, and let $( e_\alpha )_\alpha$ be a bounded approximate identity for $\A$. Note that
\begin{equation} \label{w*lim}
  f = \text{$\sigma(\A^\ast,\A)$-}\lim_\alpha e_\alpha f.
\end{equation}
By \cite[Proposition III.5.12]{Tak}, the set $\{ af : a \in \A, \, \| a \| \leq 1 \}$ is relatively weakly compact in $\A^\ast$. Consequently, $( e_\alpha )_\alpha$ has subnet $( e_\beta)_\beta$ such that $( e_\beta f )_\beta$ is weakly convergent. Since $\sigma(\A^\ast,\A^{\ast\ast})$ is finer than $\sigma(\A^\ast,\A)$, it follows from (\ref{w*lim}) that $f = \text{$\sigma(\A^\ast,\A^{\ast\ast})$-}\lim_\beta e_\beta f$. Hence, $f$ lies in the weak, i.e., norm, closure of $\{ a f : a \in \A \}$. By Cohen's factorization theorem (\cite[Corollary 2.9.26]{Dal}) there are thus $a \in \A$ and $g \in E$ such that $f = ag$.
\par 
An analogous argument applied to $g$---now with respect to the right module action of $\A$---yields $b \in \A$ and $h \in E$ such that $g = hb$.
\end{proof}
\begin{definition} \label{lucdef}
Let $\G$ be a locally compact quantum group. We define
\begin{alphitems}
\item the space of \emph{left uniformly continuous elements} of $L^\infty(\G)$ as
\[
  \LUC(\G) := \clspan\{ x \cdot f : x \in L^\infty(\G) , \, f \in L^1(\G) \},
\] 
\item the space of \emph{right uniformly continuous elements} of $L^\infty(\G)$ as
\[
  \RUC(\G) := \clspan\{ f \cdot x : x \in L^\infty(\G) , \, f \in L^1(\G) \},
\] 
and
\item the space of \emph{right uniformly continuous elements} of $L^\infty(\G)$ as 
\[
  \UC(\G) := LUC(\G) \cap \RUC(\G).
\]
\end{alphitems}
\end{definition}
\begin{remarks}
\item In view of \cite[(32.45)(a) and (b)]{HR}, these definitions are just the usual ones if $\G$ is a locally compact group. Since $A(G)$ is a commutative Banach algebra for any locally compact group $G$, we have $\LUC(\hat{G}) = \RUC(\hat{G}) = \UC(\hat{G})$; our usage of the symbol $\UC(\hat{G})$ is the same as Granirer's.
\item Let $R$ be the unitary antipode of $\G$, and note that
\[
  \begin{split}
  R(x \cdot f) & = R( ( f \tensor \id)(\Gamma x)) \\
  & = (R^\ast f \tensor \id)((R \tensor R)(\Gamma x)) \\
  & = (\id \tensor R^\ast f)(\sigma( (R \tensor R)(\Gamma x)) \sigma) \\
  & = (\id \tensor R^\ast f)(\Gamma(Rx)) = (R^\ast f) \cdot (Rx) \qquad (x \in L^\infty(\G), \, f \in L^1(\G)).
  \end{split}
\]
It follows that $R$ maps $\LUC(\G)$ isometrically onto $\RUC(\G)$ and vice versa: this fact will be useful in the sequel.
\item By (\ref{squaredense}), we have
\[
  \LUC(\G) := \clspan\{ x \cdot f : x \in \LUC(\G) , \, f \in L^1(\G) \};
\]
analogous statements hold for $\RUC(\G)$ and $\UC(\G)$.
\item The space $\RUC(\G)$ was already introduced in \cite{HNR} and is featured in \cite[Theorem 14]{HNR}. That result, however, is not connected with those obtained in this paper.
\end{remarks}
\par 
Given a locally compact group $G$, the spaces $\LUC(G)$, $\RUC(G)$, and $\UC(G)$ are obviously $\cstar$-subalgebras of ${\cal C}(G)$, and it is not difficult to see that ${\cal C}_0(G) \subset \UC(G)$. On the dual side, $\UC(\hat{G})$ is a $\cstar$-algebra containing $C^\ast_r(G)$ and contained in ${\cal M}(C^\ast_r(G))$ (\cite[Propositions 1 and 2]{Gra2}).
\par 
We shall now see that some of this carries over to general locally compact quantum groups.
\par
As in \cite{Pau}, we call a subspace of a unital $\cstar$-algebra $\A$ an \emph{operator system} in $\A$ if it contains the identity and is closed under taking adjoints. If $\A$ is a $\cstar$-algebra, and $f \in \A^\ast$, we define $\bar{f} \in \A^\ast$ by letting $\left\langle x, \bar{f} \right\rangle := \overline{\langle x^\ast, f \rangle}$ for $x \in \A$. Further, if $\xi$ and $\eta$ are vectors in a Hilbert space $\Hilbert$, we use the symbol $\omega_{\xi,\eta}$ for the functional ${\cal B}(\Hilbert) \ni T \mapsto \langle T\xi,\eta\rangle$ as well as for its restrictions to various subalgebras of ${\cal B}(\Hilbert)$; if $\eta = \xi$, we simply write $\omega_\xi$ instead of $\omega_{\xi,\xi}$.
\begin{theorem} \label{lucthm1}
Let $\G$ be a locally compact quantum group. Then $\LUC(\G)$, $\RUC(\G)$, and $\UC(\G)$ are closed operator systems in $L^\infty(\G)$ such that
\[
  {\cal C}_0(\G) \subset \UC(\G) \qquad\text{and}\qquad \LUC(\G) \cup \RUC(\G) \subset {\cal M}({\cal C}_0(\G)).
\]
\end{theorem}
\begin{proof}
Obviously, $1 \in \LUC(\G)$, and since $(x \cdot f)^\ast = x^\ast \cdot \bar{f}$ for $x \in L^\infty(\G)$ and $f \in L^1(\G)$, it is clear that $\LUC(\G)$ is an operator system in $L^\infty(\G)$. An analogous argument yields the same for $\RUC(\G)$, which entails the corresponding claim for $\UC(\G)$.
\par 
To see that ${\cal C}_0(\G) \subset \UC(\G)$, first recall (\cite[Corollary 6.11]{KV1}) that
\[
  {\cal C}_0(\G) \wTensor {\cal C}_0(\G) = \clspan \{ (\Gamma a)(1 \tensor b) : a, b \in {\cal C}_0(\G) \}.
\]
It follows that
\[
  \begin{split}
  {\cal C}_0(\G) 
  & = \{ (\id \tensor f)(\boldsymbol{a}) : f \in L^1(\G), \, \boldsymbol{a} \in {\cal C}_0(\G) \wTensor {\cal C}_0(\G) \} \\
  & = \clspan \{ (\id \tensor f)((\Gamma a)(1 \tensor b)): f \in L^1(\G), \, a, b \in {\cal C}_0(\G) \} \\
  & = \clspan \{ (\id \tensor bf)(\Gamma a): f \in L^1(\G), \, a, b \in {\cal C}_0(\G) \} \\
  & = \clspan \{ (\id \tensor f)(\Gamma a): f \in L^1(\G), \, a \in {\cal C}_0(\G) \}, \\
  & \qquad\qquad\text{by Lemma \ref{faclem} with $\A = {\cal C}_0(\G)$ and $E = L^1(\G)$}, \\
  & \subset \LUC(\G).
  \end{split}
\]
Analogously, we obtain ${\cal C}_0(\G) \subset \RUC(\G)$.
\par 
To see that $\RUC(\G) \subset {\cal M}({\cal C}_0(\G))$, let $x \in L^\infty(\G)$, $f \in L^1(\G)$, and $a \in {\cal C}_0(\G)$. We claim that $(f \cdot x) a \in {\cal C}_0(\G)$. Since $L^\infty(\G)$ is in standard form on $L^2(\G)$ (see \cite{Tak}, for instance), we may suppose that $f = \omega_{\xi,\eta}$ with $\xi,\eta \in L^2(\G)$. Let $W \in {\cal B}(L^2(\G) \ttensor_2 L^2(\G))$ be the multiplicative unitary of $\G$, and recall that $W \in {\cal M}({\cal C}_0(\G) \wTensor {\cal K}(L^2(\G)))$. Let $K \in {\cal K}(L^2(\G))$ be such that $K \xi = \xi$ (a rank one operator will do). Since $W, 1 \tensor x \in {\cal M}({\cal C}_0(\G) \wTensor {\cal K}(L^2(\G)))$, it follows that
\[
  W^\ast(1 \tensor x)W(a \tensor K) \in {\cal C}_0(\G) \wTensor {\cal K}(L^2(\G))
\]
and thus
\[
  \begin{split}
  (\omega_{\xi,\eta} \cdot x)a & = (\id \tensor \omega_{\xi,\eta})(W^\ast(1 \tensor x)W) a \\
  & = (\id \tensor \omega_{\xi,\eta})(W^\ast(1 \tensor x)W(a \tensor 1)) \\
  & = (\id \tensor \omega_{K\xi,\eta})(W^\ast(1 \tensor x)W(a \tensor 1)), \qquad\text{because $K\xi = \xi$}, \\
  & = (\id \tensor \omega_{\xi,\eta})(W^\ast(1 \tensor x)W(a \tensor K)) \in{\cal C}_0(\G).
  \end{split}
\]
It follows that $\RUC(\G) \subset {\cal RM}({\cal C}_0(\G))$; applying the involution, yields $\RUC(\G) \subset {\cal LM}({\cal C}_0(\G))$ as well and thus $\RUC(\G) \subset {\cal M}({\cal C}_0(\G))$.
\par
To show that $\LUC(\G) \subset {\cal M}({\cal C}_0(\G))$ as well, just note that the unitary antipode $R$ of $\G$ leaves ${\cal C}_0(\G)$ invariant, reverses multiplication---and thus leaves ${\cal M}({\cal C}_0(\G))$ invariant---, and maps $\LUC(\G)$ bijectively onto $\RUC(\G)$, as was remarked immediately after Definition \ref{lucdef}.
\end{proof}
\begin{remark}
Unless $\G$ is discrete, i.e., $L^1(\G)$ has an identity, we cannot suppose that $L^1(\G) = M(\G)$. Hence, we cannot use the density condition of \cite[Definition 4.1]{KV1} directly to conclude that ${\cal C}_0(\G) \subset \UC(\G)$. 
\end{remark}
\par 
Theorem \ref{lucthm1} leaves open the question of whether $\LUC(\G)$, $\RUC(\G)$, and $\UC(\G)$ are, in fact, $\cstar$-subalgebras of $L^\infty(\G)$ as opposed to mere operator systems. We shall deal with this question in Section \ref{cstar?} below.
\section{Amenability and co-amenability}
The notions of invariant means and amenability are well established for locally compact groups (see \cite{Pie}). These notions extend naturally to locally compact quantum groups:
\begin{definition} \label{LIMdef}
Let $\G$ be a locally compact quantum group, and let $E$ be an operator system in $L^\infty(\G)$ which is also a right $L^1(\G)$-submodule of $L^\infty(\G)$. Then $M \in E^\ast$ with $\| M \| = \langle 1, M \rangle = 1$ is called a \emph{left invariant mean} on $E$ if
\[
  f \cdot M = \langle 1, f \rangle M \qquad (f \in L^1(\G)).
\]
If there is a left invariant mean on $L^\infty(\G)$, we call $\G$ \emph{amenable}.
\end{definition}
\begin{remarks}
\item The condition that $\| M \| = \langle 1, M \rangle = 1$ is equivalent to $M$ being positive with $\langle 1, M \rangle = 1$ (\cite[Proposition 2.11 and Exercise 2.3]{Pau}).
\item Our choice of terminology is in accordance with \cite{BT}, but not universally agreed upon: amenable, locally compact quantum groups---or, rather, Kac algebras---in the sense of Definition \ref{LIMdef} are called \emph{Voiculescu amenable} in \cite{Rua} and \emph{weakly amenable} in \cite{DQV}.
\item In view of Definition \ref{LIMdef}, it is obvious what the definitions of a right invariant or simply (two-sided) invariant mean are supposed to be. The existence of a right invariant and of an invariant mean on $L^\infty(\G)$ is equivalent to the amenability of $\G$ (see \cite{DQV}).
\end{remarks}
\par 
For a locally compact group $G$, the amenability of $G$ is equivalent to the existence of various types of invariant means on various subspaces of $L^\infty(G)$ (see again \cite{Pie}); in particular, $G$ is amenable if and only if there is a left invariant mean---in the sense of Definition \ref{LIMdef}, which would be called a topologically invariant mean in \cite{Pie}---on ${\cal C}(G)$, $\LUC(G)$, $\RUC(G)$, or $\UC(G)$.
\par 
In \cite{BT}, B\'edos and Tuset consider the existence of a left invariant mean on ${\cal M}({\cal C}_0(\G))$ (a condition they term, somewhat misleadingly, \emph{topological amenability}), and observe that this condition is formally weaker than amenability. They remark that, unless in the trivial case when $\G$ is discrete, it is not clear whether the existence of such a left invariant mean is equivalent to amenability, as it is in the group case.
\par 
Making use of Theorem \ref{lucthm1}, we shall see that, for a considerably larger class of locally compact quantum groups, the existence of a left invariant mean on ${\cal M}({\cal C}_0(\G))$ does indeed imply amenability.
\begin{definition}
A locally compact quantum group $\G$ is called \emph{co-amenable} if $L^1(\G)$ has a bounded approximate identity.
\end{definition}
\begin{remarks}
\item Trivially, every discrete quantum group is co-amenable, but so is every locally compact group $G$ whereas $\hat{G}$ is co-amenable if and only if $G$ is amenable (\cite{Lep}).
\item If $\G$ is co-amenable, then $\hat{\G}$ is amenable (\cite{BT}). Whether the converse holds is a major open problem: it is known to be true for groups (\cite{Lep}) and discrete quantum groups (\cite{Tom}).
\item If $\G$ is co-amenable, then Cohen's factorization theorem (\cite[Corollary 2.9.26]{Dal}) yields that
\[
  \LUC(\G) = \{ x \cdot f : f \in L^1(\G), \, x \in L^\infty(\G) \};
\]
analogous statements hold for $\RUC(\G)$ and $\UC(\G)$.
\end{remarks}
\par 
By \cite[Theorem 3.1]{BT}, the existence of a one-sided approximate identity for $L^1(\G)$ is already enough to ascertain the co-amenability of $\G$. An inspection of the proof of that theorem shows that such one-sided approximate identities can be chosen to consist of states. Somewhat less obvious is the fact that, if $\G$ is co-amenable, a two-sided approximate identity for $L^1(\G)$ can be found that consists of states (\cite[Theorem 2]{HNR}). Even though it is not directly related to our investigation of uniform continuity, we note the following improvement of \cite[Theorem 2]{HNR}, which is surprising even in the group case:
\begin{theorem} \label{lucthm2}
Let $\G$ be a locally compact quantum group. Then the following are equivalent for a net $( \xi_\alpha )_\alpha$ of unit vectors in $L^2(\G)$:
\begin{items}
\item $(\omega_{\xi_\alpha} )_\alpha$ is a bounded left approximate identity for $L^1(\G)$;
\item $\| W(\xi_\alpha \tensor \eta) - \xi_\alpha \tensor \eta \| \to 0 \qquad (\eta \in L^2(\G))$;
\item $(\omega_{\xi_\alpha} )_\alpha$ is a bounded (two-sided) approximate identity for $L^1(\G)$.
\end{items}
In particular, every left approximate identity for $L^1(\G)$ consisting of states is a bounded approximate identity for $L^1(\G)$.
\end{theorem}
\begin{proof}
(i) $\Longleftrightarrow$ (ii) is part of the proof of \cite[Theorem 3.1]{BT}, and (iii) $\Longrightarrow$ (i) is trivial.
\par 
(ii) $\Longrightarrow$ (iii): Suppose that (ii) holds. By (ii) $\Longrightarrow$ (i), we already know that $( \omega_{\xi_\alpha} )_\alpha$ is a bounded left approximate identity for $L^1(\G)$. Hence, it is enough to show that $( \omega_{\xi_\alpha} )_\alpha$ is a bounded right approximate identity for $L^1(\G)$. 
\par 
Let $\G^\op$ denote the opposite quantum group of $\G$ as defined in \cite[Definition 4.1]{KV2}. By \cite[p.\ 90]{KV2} and \cite[(2.6)]{JNR}, 
\[
  W^\op = \sigma V^\ast \sigma = (1 \tensor \upsilon) W^\ast (1 \tensor \upsilon^\ast)
\]
holds with $\upsilon := \hat{J}J$, where $\hat{J}$ and $J$ are the modular conjugations obtained from $\phi$ and $\hat{\phi}$, respectively. It thus follows from (ii) that
\begin{multline*}
  \| W^\op(\xi_\alpha \tensor \eta) - \xi_\alpha \tensor \eta \| \\ = 
  \| (1 \tensor \upsilon) W^\ast(\xi_\alpha \tensor \upsilon^\ast \eta) - \xi_\alpha \tensor \eta \| =
  \| \xi_\alpha \tensor \upsilon^\ast \eta - W(\xi_\alpha \tensor \upsilon^\ast \eta) \| \to 0 \qquad (\eta \in L^2(\G)).
\end{multline*}
From (ii) $\Longrightarrow$ (i)---applied to $\G^\op$ instead of $\G$---, we conclude that $( \omega_{\xi_\alpha} )_\alpha$ is a bounded left approximate identity for $L^1(\G^\op)$. From the definition of $\G^\op$, it is obvious that $L^1(\G^\op)$ is nothing but $L^1(\G)$ with reversed multiplication. Hence, $( \omega_{\xi_\alpha} )_\alpha$ is a bounded right approximate identity for $L^1(\G)$.
\end{proof}
\par 
We now present our main result in this section:
\begin{theorem} \label{lucthm3}
Let $\G$ be a co-amenable, locally compact quantum group. Then the following are equivalent:
\begin{items}
\item $\G$ is amenable;
\item there is a left invariant mean on ${\cal M}({\cal C}_0(\G))$;
\item there is a left invariant mean on $\LUC(\G)$;
\item there is a left invariant mean on $\RUC(\G)$;
\item there is a left invariant mean on $\UC(\G)$.
\end{items}
\end{theorem}
\begin{proof}
(i) $\Longrightarrow$ (ii) $\Longrightarrow$ (iii) $\Longrightarrow$ (v) and (ii) $\Longrightarrow$ (iv) $\Longrightarrow$ (v) are obvious in the light of Theorem \ref{lucthm1}: just restrict the respective left invariant mean to the smaller space.
\par 
(v) $\Longrightarrow$ (i): Let $\tilde{M}$ be a left invariant mean on $\UC(\G)$, and let $( e_\alpha )_\alpha$ be an approximate identity for $L^1(\G)$ consisting of states. Let $\cal U$ be an ultrafilter on the index set of $( e_\alpha )_\alpha$ that dominates the order filter. Define $M \!: L^\infty(\G) \to \comps$ by letting
\[
  \langle x, M \rangle := 
  \lim_{\alpha \in \cal U} \left\langle e_\alpha \cdot x \cdot e_\alpha, \tilde{M} \right\rangle \qquad (x \in L^\infty(\G)).
\]
It is immediate that $\| M \| \leq 1$ and $\langle 1, M \rangle = \left\langle 1, \tilde{M} \right\rangle = 1$, so that $M$ is a state on $L^\infty(\G)$. Moreover, note that
\[
  \begin{split}
  \langle x, f \cdot M \rangle & = \langle x \cdot f, M \rangle \\
  & = \lim_{\alpha \in \cal U} \left\langle e_\alpha \cdot (x \cdot f) \cdot e_\alpha, \tilde{M} \right\rangle \\
  & = \lim_{\alpha \in \cal U} \left\langle e_\alpha \cdot x \cdot (f \ast e_\alpha), \tilde{M} \right\rangle \\
  & = \lim_{\alpha \in \cal U} \left\langle e_\alpha \cdot x \cdot (e_\alpha \ast f), \tilde{M} \right\rangle \\
  & = \lim_{\alpha \in \cal U} \left\langle (e_\alpha \cdot x \cdot e_\alpha) \cdot f, \tilde{M} \right\rangle \\
  & = \lim_{\alpha \in \cal U} \left\langle e_\alpha \cdot x \cdot e_\alpha, f \cdot \tilde{M} \right\rangle \\
  & = \lim_{\alpha \in \cal U} \langle 1, f \rangle \left\langle e_\alpha \cdot x \cdot e_\alpha, \tilde{M} \right\rangle \\
  & = \langle 1, f \rangle \langle x, M \rangle,
  \end{split}
\]
so that $M$ is a left invariant mean on $L^\infty(\G)$.
\end{proof}
\begin{remarks}
\item It is easy to state (and prove) a right or two-sided version of Theorem \ref{lucthm3}.
\item The proof of (v) $\Longrightarrow$ (i) in Theorem \ref{lucthm3} does require the existence of a net $( e_\alpha )_\alpha$ of states in $L^1(\G)$ such that
\[
  f \ast e_\alpha - e_\alpha \ast f \to 0 \qquad (f \in L^1(\G)),
\]
which is a condition weaker than co-amenability. It is also satisfied by every ame\-na\-ble, locally compact quantum group---a somewhat pointless hypothesis for Theorem \ref{lucthm3}---and may be true for every locally compact quantum group.
\end{remarks}
\section{Uniform continuity and weak almost periodicity}
Recall that a bounded, continuous function on a locally compact group is called \emph{weakly almost periodic} if the set $\{ L_x f : x \in G \}$ is relatively weakly compact in ${\cal C}(G)$. Every function in ${\cal C}_0(G)$ is weakly almost periodic (\cite[Corollary 3.7]{Bur}) and every weakly almost periodic function is uniformly continuous (\cite[Theorem 3.11]{Bur}). More generally, one can consider the space of weakly almost periodic functionals on any Banach algebra (see \cite{Lau1}, \cite{LL}, or \cite{Run2}, for instance): the weakly almost periodic function on $G$ then correspond to the weakly almost periodic functionals on $L^1(G)$.
\par 
Specializing the Banach algebraic definition to $L^1(\G)$ for a locally compact quantum group $\G$, we define:
\begin{definition} \label{wapdef}
Let $\G$ be a locally compact quantum group. We define the space of \emph{weakly almost periodic elements} of $L^\infty(\G)$ as 
\[
  \WAP(\G) := \{ x \in L^\infty(\G) : \text{$L^1(\G) \ni f \mapsto f \cdot x$ is weakly compact} \}.
\]
\end{definition}
\begin{remarks}
\item For a locally compact group $G$, the space $\WAP(G)$ in the sense of Definition \ref{wapdef} just consists of the weakly almost periodic functions on $G$ studied in \cite{Bur}.
\item Even though Definition \ref{wapdef} appears to be somewhat asymmetric, Grothendieck's double limit criterion immediately yields that
\[
  \WAP(\G) = \{ x \in L^\infty(\G) : \text{$L^1(\G) \ni f \mapsto x \cdot f$ is weakly compact} \}
\]
as well.
\end{remarks}
\par 
In \cite{Run1}, a Banach algebra $\A$ was called \emph{dual} if there is a---not necessarily unique---Banach space $\A_\ast$ with $\A = (\A_\ast)^\ast$ such that multiplication in $\A$ is separately $\sigma(\A,\A_\ast)$-continuous. For instance, if $G$ is a locally compact group, then $M(G) = {\cal C}_0(\G)^\ast$ is a dual Banach algebras. This extends to locally compact quantum groups:
\begin{proposition} \label{dualprop}
Let $\G$ be a locally compact quantum group. Then $M(\G) = {\cal C}_0(\G)^\ast$ is a dual Banach algebra.
\end{proposition}
\begin{proof}
Let $\nu \in M(\G)$, and let $a \in {\cal C}_0(\G)$.
\par 
We claim that
\begin{equation} \label{w*cont}
  M(\G) \to \comps, \quad \mu \mapsto \langle a, \mu \ast \nu \rangle 
\end{equation}
is weak$^\ast$ continuous. By Lemma \ref{faclem}, there are $b \in {\cal C}_0(\G)$ and $\tilde{\nu} \in M(\G)$ such that $\nu = b\tilde{\nu}$. Note that
\begin{equation} \label{slice}
  \langle a, \mu \ast \nu \rangle = \langle \Gamma a, \mu \tensor \nu \rangle = 
  \left\langle \Gamma a, \mu \tensor b \tilde{\nu} \right\rangle
  = \left\langle \left(\id \tensor \tilde{\nu}\right) ((\Gamma a)(1 \tensor b)), \mu \right\rangle.
\end{equation}
By \cite[Corollary 6.11]{KV1}, we have $(\Gamma a)(1 \tensor b) \in {\cal C}_0(\G) \wTensor {\cal C}_0(\G)$ and thus $\left(\id \tensor \tilde{\nu}\right) ((\Gamma a)(1 \tensor b)) \in {\cal C}_0(\G)$. In view of (\ref{slice}), this yields the weak$^\ast$ continuity of (\ref{w*cont}). 
\par
Analogously, we see that
\[
  M(\G) \to \comps, \quad \mu \mapsto \langle a, \nu \ast \mu \rangle 
\]
is weak$^\ast$ continuous.
\end{proof}
\par
Simultaneously extending \cite[Corollary 3.7]{Bur} on the one hand and \cite[Theorem 3.11]{Bur} and \cite[Proposition 1]{Gra1} on the other, we obtain:
\begin{theorem} \label{wapthm}
Let $\G$ be locally compact quantum group. Then $\WAP(\G)$ is a closed operator system in $L^\infty(\G)$ containing ${\cal C}_0(\G)$. Moreover, if $\G$ is co-amenable, then $\WAP(\G) \subset \UC(\G)$ holds.
\end{theorem}
\begin{proof}
It is straightforward to verify that $\WAP(\G)$ is a closed operator system in $L^\infty(\G)$.
\par 
Let $a \in {\cal C}_0(\G)$. By Proposition \ref{dualprop}, the map
\begin{equation} \label{C0}
   M(\G) \to {\cal C}_0(\G), \quad \mu \mapsto \mu \cdot a 
\end{equation}
is weak$^\ast$-weakly continuous. Since the closed unit ball of $M(\G)$ is weak$^\ast$ compact, this means that (\ref{C0}) is weakly compact, as is its restriction to $L^1(\G)$.
\par
Suppose that $\G$ is co-amenable, and let $x \in \WAP(\G)$. An argument almost identical to that in the proof of Lemma \ref{faclem} reveals that $x$ lies in the norm closure of $\{ f \cdot x : f \in L^1(\G) \}$ and thus in $\RUC(\G)$. Analogously, one shows that $\WAP(\G) \subset \LUC(\G)$ as well.
\end{proof}
\begin{remark}
Together, Theorems \ref{lucthm1} and \ref{wapthm} yield for co-amenable $\G$ that
\begin{equation} \label{lucwapincl}
  {\cal C}_0(\G) \subset \WAP(\G) \subset \UC(\G) \subset \LUC(\G) \cup \RUC(\G) \subset {\cal M}({\cal C}_0(\G)).
\end{equation}
If, in addition, $\G$ is compact, i.e., ${\cal C}_0(\G)$ has an identity, we have ${\cal C}_0(\G) = {\cal M}({\cal C}_0(\G))$, and thus equality throughout in (\ref{lucwapincl}); in particular, $\LUC(\G) = \WAP(\G)$ holds. For a locally compact group $G$, the equality $\LUC(G) = \WAP(G)$ is, in fact, equivalent to $G$ being compact: this follows from A.\ T.-M.\ Lau's description of the topological center of $\LUC(G)^\ast$ (\cite{Lau2}). Whether this equivalence is also true for general---or at least for co-amenable---locally compact quantum group appears to be wide open.
\end{remark}
\par 
If $G$ is a locally compact group, then there is a unique invariant mean on $\WAP(G)$ (\cite[Theorem 1.25]{Bur}), and the same is true for $\WAP(\hat{G})$ (\cite[Corollary 5.7]{Lau1}). For amenable, locally compact quantum groups, we obtain:
\begin{proposition} \label{wapprop}
Let $\G$ be an amenable, locally compact quantum group. Then there is a unique left invariant mean on $\WAP(\G)$ which is automatically right invariant. 
\end{proposition}
\begin{proof}
By \cite[Theorem 4.10]{Run2}, $\WAP(\G)^\ast$ is a dual Banach algebra in a canonical fashion; we denote its product by $\ast$ (restricted to $L^1(\G)$, it is just the product defined in (\ref{prod})). It is routine to check that $M \in \WAP(\G)^\ast$ with $\| M \| = \langle 1, M \rangle = 1$ is a left invariant mean on $\WAP(\G)$ if and only if
\begin{equation} \label{LIM}
  F \ast M = \langle 1,F \rangle M \qquad (F \in \WAP(\G)^\ast)
\end{equation}
and a right invariant mean if and only if
\begin{equation} \label{RIM}
  M \ast F = \langle 1,F \rangle M \qquad (F \in \WAP(\G)^\ast).
\end{equation}
\par 
As $\G$ is amenable, there is a (two-sided) invariant mean on $L^\infty(\G)$ the restriction of which to $\WAP(\G)$ is an invariant mean, say $M_0$, on $\WAP(\G)$. Let $M$ be any left invariant mean on $\WAP(\G)$. Then (\ref{LIM}) and (\ref{RIM}) yield
\[
  M_0 = \langle 1, M \rangle M_0 = M_0 \ast M = \langle 1, M_0 \rangle M = M,
\]
which completes the proof.
\end{proof}
\begin{remark}
There is an invariant mean on $\WAP(G)$ for \emph{any} locally compact group $G$, i.e., without any amenability hypothesis. It is an interesting question whether the same is true for locally compact quantum groups. The proof of Proposition \ref{wapprop} shows that whenever there is a left invariant mean on $\WAP(\G)$, it is necessarily unique and also right invariant.
\end{remark}
\section{$\LUC(\G)$ and $\RUC(\G)$ as $\cstar$-algebras} \label{cstar?}
If $G$ is a locally compact group, then $\LUC(G)$, $\RUC(G)$, and $\UC(G)$ are obviously $\cstar$-algebras. Somewhat less obvious is that $\UC(\hat{G})$ is also a $\cstar$-algebra (see \cite{Gra2}). For a locally compact quantum group $\G$ it is easy to see that $\LUC(\G)$, $\RUC(\G)$, and $\UC(\G)$ are $\cstar$-algebras if $\G$ is discrete or compact: in the discrete case, $\LUC(\G) = \RUC(\G) = L^\infty(\G)$ holds trivially, and if $\G$ is compact, Theorem \ref{lucthm1} yields $\LUC(\G) = \RUC(\G) = {\cal C}_0(\G)$.
\par 
In this section, we give, for co-amenable $\G$, an alternative description of $\LUC(\G)$ and $\RUC(\G)$ as spaces of quasi-multipliers. Under another technical condition, we shall see that $\LUC(\G)$ and $\RUC(\G)$ consist even of multipliers and are indeed $\cstar$-algebras. 
\par 
We begin with a lemma, which we formulate using leg notation. To make its proof less cumbersome to formulate, we also introduce some notation related to operator spaces: given any two operator spaces $E$ and $F$, we denote by ${\cal CA}(E,F)$ the closure of the finite rank operators in $\CB(E,F)$, which can be canonically identified with the injective tensor product $F \wTensor E^\ast$ (\cite[Proposition 8.1.2]{ER}). 
\begin{lemma} \label{ruclem}
Let $\A$ be a $\cstar$-algebra, let $\Hilbert$ be a Hilbert space, let $A, B \in \A \wTensor {\cal K}(\Hilbert)$, let $T \in {\cal B}(\Hilbert) \bar{\tensor} {\cal B}(\Hilbert)$, and let $\xi, \eta \in \Hilbert$. Then we have
\begin{equation} \label{CA0}
  (\id \tensor \id \tensor \omega_{\xi,\eta})(A_{2,3} T_{1,3} B_{2,3}) \in {\cal B}(\Hilbert) \wTensor \A.
\end{equation}
\end{lemma}
\begin{proof}
By \cite[Lemma 3.1]{DR}, the map
\[
  {\cal B}(\Hilbert) \to \A \wTensor {\cal K}(\Hilbert), \quad x \mapsto A(1 \tensor x)B
\]
lies in ${\cal CA}({\cal B}(\Hilbert),\A \wTensor {\cal K}(\Hilbert))$, so that
\begin{equation} \label{CA1}
  {\cal B}(\Hilbert) \to \A, \quad x \mapsto (\id \tensor \omega_{\xi,\eta})(A(1 \tensor x)B)
\end{equation}
belongs to ${\cal CA}({\cal B}(\Hilbert),\A)$.
\par 
The canonical completely isometric isomorphisms
\[
  {\cal B}(\Hilbert) \bar{\tensor} {\cal B}(\Hilbert) \cong ( {\cal B}(\Hilbert)_\ast \Tensor {\cal B}(\Hilbert)_\ast )^\ast
  \cong \CB({\cal B}(\Hilbert)_\ast,{\cal B}(\Hilbert)) 
\]
by \cite[Corollary 7.1.5 and Theorem 7.2.4]{ER} allow us to interpret $T$ as an element of $\CB({\cal B}(\Hilbert)_\ast,{\cal B}(\Hilbert))$ via
\begin{equation} \label{CA2}
  {\cal B}(\Hilbert)_\ast \to {\cal B}(\Hilbert), \quad \omega \mapsto (\id \tensor \omega)(T).
\end{equation}
Composing (\ref{CA2}) with (\ref{CA1}) then yields an operator in ${\cal CA}({\cal B}(\Hilbert)_\ast,\A)$, namely
\[
  {\cal B}(\Hilbert)_\ast \to \A, \quad
  \omega \mapsto (\omega \tensor \id \tensor \omega_{\xi,\eta})(A_{2,3} T_{1,3} B_{2,3}).
\]
That this operator lies in ${\cal CA}({\cal B}(\Hilbert)_\ast,\A)$ is obviously equivalent to (\ref{CA0}).
\end{proof}
\begin{proposition} \label{lucprop}
Let $\G$ be a locally compact quantum group. Then 
\begin{equation} \label{lucincl}
  \LUC(\G) \subset \{ x \in {\cal M}({\cal C}_0(\G)) : 
  \text{$(a \tensor 1)(\Gamma x)(b \tensor 1) \in {\cal C}_0(\G) \wTensor {\cal B}(L^2(\G))$ for all $a,b \in {\cal C}_0(\G)$} \} 
\end{equation}
and
\begin{equation} \label{rucincl}
  \RUC(\G) \subset \{ x \in {\cal M}({\cal C}_0(\G)) : 
  \text{$(1 \tensor a)(\Gamma x)(1 \tensor b) \in {\cal B}(L^2(\G)) \wTensor {\cal C}_0(\G)$ for all $a,b \in {\cal C}_0(\G)$} \} 
\end{equation}
hold, with equality in both cases if $\G$ is co-amenable.
\end{proposition}
\begin{proof}
We shall prove (\ref{rucincl}); applying the unitary antipode---or, alternatively, simply considering $\G^\op$---, then yields (\ref{lucincl}).
\par
Let $x \in L^\infty(\G)$, and let $f \in L^1(\G)$; note that
\begin{equation} \label{ruceq1}
  \begin{split}
  \Gamma(f \cdot x) & = \Gamma( (\id \tensor f)(\Gamma x)) \\
  & = (\Gamma \tensor f)(\Gamma x) \\
  & = (\id \tensor \id \tensor f)((\Gamma \tensor \id)(\Gamma x)) \\
  & = (\id \tensor \id \tensor f)((\id \tensor \Gamma)(\Gamma x)) \\
  & = (\id \tensor \id \tensor f)(W_{2,3}^\ast W_{1,3}^\ast (1 \tensor 1 \tensor x) W_{1,3} W_{2,3}).
  \end{split}
\end{equation}
Let $a,b \in {\cal C}_0(\G)$. Choose $\xi, \eta \in L^2(\G)$ such that $f = \omega_{\xi,\eta}$, which can be done because $L^\infty(\G)$ is in standard form on $L^2(\G)$, and furthermore, pick $K,L \in {\cal K}(L^2(\G))$ such that $K\xi = \xi$ and $L^\ast \eta = \eta$. In view of (\ref{ruceq1}), we obtain
\[
  \begin{split}
  \lefteqn{(1 \tensor a)(\Gamma(f \cdot x))(1 \tensor b)} \\
  & = (1 \tensor a)((\id \tensor \id \tensor \omega_{\xi,\eta})
  (W_{2,3}^\ast W_{1,3}^\ast (1 \tensor 1 \tensor x) W_{1,3} W_{2,3}))(1 \tensor b) \\
  & = (\id \tensor \id \tensor \omega_{K\xi,L^\ast\eta}) 
  ((1 \tensor a \tensor 1)(W_{2,3}^\ast W_{1,3}^\ast (1 \tensor 1 \tensor x) W_{1,3} W_{2,3})(1 \tensor b \tensor 1)) \\
  & = (\id \tensor \id \tensor \omega_{\xi,\eta})
  ((1 \tensor a \tensor L)(W_{2,3}^\ast W_{1,3}^\ast (1 \tensor 1 \tensor x) W_{1,3} W_{2,3})(1 \tensor b \tensor K)) \\
  & = (\id \tensor \id \tensor \omega_{\xi,\eta})(A_{2,3} T_{1,3} B_{2,3}),
  \end{split}
\]
where
\begin{align*}
  A & := (a \tensor L)W^\ast \in {\cal C}_0(\G) \wTensor {\cal K}(L^2(\G)), \\
  B & := W(b \tensor K) \in {\cal C}_0(\G) \wTensor {\cal K}(L^2(\G)), \\
\intertext{and}
  T & := W^\ast (1 \tensor x) W \in {\cal B}(L^2(\G)) \bar{\tensor} {\cal B}(L^2(\G)).
\end{align*}
From Lemma \ref{ruclem}, we conclude that $(1 \tensor a)(\Gamma(f \cdot x))(1 \tensor b) \in {\cal B}(L^2(\G)) \wTensor {\cal C}_0(\G)$. In view of the definition of $\RUC(\G)$ and Theorem \ref{lucthm1}, we obtain (\ref{rucincl}).
\par 
Suppose now that $\G$ is co-amenable. We shall prove that the inclusion (\ref{lucincl}) is, in fact, an equality (obtaining the equality in (\ref{rucincl}) is then easy again).
\par 
Let $x \in {\cal M}({\cal C}_0(\G))$ be such that $(a \tensor 1)(\Gamma x)(b \tensor 1) \in {\cal C}_0(\G) \wTensor {\cal B}(L^2(\G))$ for all $a,b \in {\cal C}_0(\G)$. Let $( e_\alpha )_{\alpha \mathbb A}$ be a bounded approximate identity for $L^1(\G)$, and let $\epsilon \in M(\G)$ be a $\sigma(M(\G),{\cal C}_0(\G))$ accumulation point of $( e_\alpha )_{\alpha \in \mathbb A}$; without loss of generality, we suppose that $\epsilon = \text{$\sigma(M(\G),{\cal C}_0(\G))$-}\lim_\alpha e_\alpha$. By (the proof of) \cite[Theorem 3.1]{BT}, $\epsilon$ is the identity element of $M(\G)$, i.e.,
\begin{equation} \label{sliceid}
  (\epsilon \tensor \id)(\Gamma a) = a \qquad (a \in {\cal C}_0(\G));
\end{equation}
furthermore, $\epsilon$ is a character of ${\cal C}_0(\G)$ (and thus on ${\cal M}({\cal C}_0(\G))$, too). By strict continuity, (\ref{sliceid}) holds not only for $a \in {\cal C}_0(\G)$, but for all $a \in {\cal M}({\cal C}_0(\G))$ as well; in particular, we have $(\epsilon \tensor \id)(\Gamma x) = x$. Let $a \in {\cal C}_0(\G)$ be such that $\langle a, \epsilon \rangle = 1$. Since $\epsilon$ is multiplicative, this means that
\[
  x = ((\epsilon \tensor \id)(a \tensor 1))((\epsilon \tensor \id)(\Gamma x))((\epsilon \tensor \id)(a \tensor 1)) 
  = (\epsilon \tensor \id)((a \tensor 1)(\Gamma x)(a \tensor 1)).
\]
By the hypothesis on $x$, we have $(a \tensor 1)(\Gamma x)(a \tensor 1) \in {\cal C}_0(\G) \wTensor {\cal B}(L^2(\G))$, so that
\begin{equation} \label{wastlim}
  x = \text{$\sigma(M(\G),{\cal C}_0(\G))$-}\lim_\alpha (e_\alpha \tensor \id)((a \tensor 1)(\Gamma x)(a \tensor 1)) 
\end{equation}
We may consider $(a \tensor 1)(\Gamma x)(a \tensor 1)$ as an element of ${\cal CA}(M(\G),{\cal B}(L^2(\G)))$ by \cite[Proposition 8.1.1]{ER}. A moment's thought reveals that those operators in ${\cal CA}(M(\G),{\cal B}(L^2(\G)))$ that arise from elements of ${\cal C}_0(\G) \wTensor {\cal B}(L^2(\G))$ have to be $\sigma(M(\G),{\cal C}_0(\G))$-norm continuous. Hence, the limit in (\ref{wastlim}) is, in fact, a norm limit. Since
\[
  (e_\alpha \tensor \id)((a \tensor 1)(\Gamma x)(a \tensor 1)) = (a e_\alpha a \tensor \id)(\Gamma x) =
  x \cdot (a e_\alpha a ) \in \LUC(\G) \qquad (\alpha \in \mathbb{A}) 
\]
this means that $x \in \LUC(\G)$.
\end{proof}
\par 
Since ${\cal B}(L^2(\G))$ has an identity, it is immediate from Proposition \ref{lucprop} that, for co-amenable $\G$, 
\[
  \LUC(\G) = \{ x \in {\cal M}({\cal C}_0(\G)) : \Gamma x \in {\cal QM}({\cal C}_0(\G) \wTensor {\cal B}(L^2(\G))) \}.
\]
and
\[
  \RUC(\G) = \{ x \in {\cal M}({\cal C}_0(\G)) : \Gamma x \in {\cal QM}({\cal B}(L^2(\G)) \wTensor{\cal C}_0(\G)) \}.
\]
As ${\cal QM}({\cal C}_0(\G) \wTensor {\cal B}(L^2(\G)))$ need not be a $\cstar$-algebra, this is not enough to conclude that $\LUC(\G)$ and $\RUC(\G)$ are $\cstar$-algebras. 
\par 
With an additional hypothesis, however, we obtain:
\begin{theorem} \label{lucthm4}
Let $\G$ be a co-amenable, locally compact quantum group, and suppose that ${\cal C}_0(\G)$ has a bounded approximate identity in its center. Then
\[
  \LUC(\G) = \{ x \in {\cal M}({\cal C}_0(\G)) : \Gamma x \in {\cal M}({\cal C}_0(\G) \wTensor {\cal B}(L^2(\G))) \}
\]
and
\[
  \RUC(\G) = \{ x \in {\cal M}({\cal C}_0(\G)) : \Gamma x \in {\cal M}({\cal B}(L^2(\G)) \wTensor{\cal C}_0(\G)) \}.
\]
holds; in particular, $\LUC(\G)$ and $\RUC(\G)$ are $\cstar$-subalgebras of ${\cal M}({\cal C}_0(\G))$.
\end{theorem}
\begin{proof}
We only prove the claim for $\LUC(\G)$.
\par
In view of the comments following Proposition \ref{lucprop}, any $x \in {\cal M}({\cal C}_0(\G))$ with $\Gamma x \in {\cal M}({\cal C}_0(\G) \wTensor {\cal B}(L^2(\G)))$ lies in $\LUC(\G)$.
\par 
Conversely, let $x \in \LUC(\G)$. 
\par 
We claim that $(\Gamma x)(a \tensor 1) \in {\cal C}_0(\G) \wTensor {\cal B}(L^2(\G))$ for each $a \in {\cal C}_0(\G)$. Let $a \in {\cal C}_0(\G)$, and let $( e_\alpha )_{\alpha \in \mathbb A}$ be a bounded approximate identity in the center of ${\cal C}_0(\G)$. Clearly, the net $(e _\alpha \tensor 1)_{\alpha \in \mathbb A}$ lies in ${\cal M}({\cal C}_0(\G) \wTensor {\cal K}(L^2(\G)))$ and commutes with every element of ${\cal C}_0(\G) \wTensor {\cal K}(L^2(\G))$; consequently, it lies in the center of ${\cal M}({\cal C}_0(\G) \wTensor {\cal K}(L^2(\G)))$. Since $(e_\alpha \tensor 1)(\Gamma x)(a \tensor 1) \in {\cal C}_0(\G) \wTensor {\cal B}(L^2(\G))$ for each $\alpha \in \mathbb A$ by Proposition \ref{lucprop}, and since $W \in {\cal M}({\cal C}_0(\G) \wTensor {\cal K}(L^2(\G)))$, we have
\[
  \begin{split}
  (\Gamma x) (a \tensor 1) & = \lim_\alpha (\Gamma x)(e_\alpha a \tensor 1) \\ 
  & = \lim_\alpha W^\ast(1 \tensor x)W(e_\alpha \tensor 1)(a \tensor 1) \\
  & = \lim_\alpha (e_\alpha \tensor 1)W^\ast(1 \tensor x)W(a \tensor 1) \\
  & = \lim_\alpha (e_\alpha \tensor 1)(\Gamma x)(a \tensor 1) \in {\cal C}_0(\G) \wTensor {\cal B}(L^2(\G)),
  \end{split}
\]
as claimed.
\par
Since $x^\ast \in \LUC(\G)$ as well, the previous argument---now applied to $x^\ast$---yields that $(\Gamma x^\ast)(a \tensor 1)$ also lies in ${\cal C}_0(\G) \wTensor {\cal B}(L^2(\G))$ for each $a \in {\cal C}_0(\G)$. Applying the involution, we conclude that $(a \tensor 1) (\Gamma x) \in {\cal C}_0(\G) \wTensor {\cal B}(L^2(\G))$ for each $a \in {\cal C}_0(\G)$ as well. Since ${\cal B}(L^2(\G))$ is unital, this is enough to ensure that $\Gamma x \in {\cal M}({\cal C}_0(\G) \wTensor {\cal B}(L^2(\G)))$.
\end{proof}
\begin{remarks}
\item If $G$ is a locally compact group, then ${\cal C}_0(G)$ trivially has a bounded approximate identity in its center whereas this is true for ${\cal C}_0(\hat{G}) = C^\ast_r(G)$ if and only if $G$ has small invariant neighborhoods (see \cite{Los}). As $\UC(\hat{G})$ is a $\cstar$-algebra for every $G$, this suggests that Theorem \ref{lucthm4} is not optimal. 
\item The hypothesis of Theorem \ref{lucthm4} that ${\cal C}_0(\G)$ have a bounded approximate identity in its center can be weakened: an inspection of the proof shows that it is sufficient for ${\cal C}_0(\G)$ to have a bounded approximate identity $( e_\alpha )_\alpha$ such that
\[
  (e_\alpha \tensor 1)W - W(e_\alpha \tensor 1) \to 0.
\]
This hypothesis may be satisfied by every (co-amenable) locally compact quantum group.
\item At the beginning of this section, we remarked that, for discrete or compact $\G$, we have $\LUC(\G) = {\cal M}({\cal C}_0(\G))$. For a locally compact group $G$, the equality $\LUC(G) = {\cal C}(G) = {\cal M}({\cal C}_0(G))$ is not only necessary, but also sufficient for $G$ to be discrete or compact: this follows from \cite[Corollary 4]{BB} and the discussion following it. It is an intriguing (and probably very hard) question whether $\LUC(\G) = {\cal M}({\cal C}_0(\G))$ forces a locally compact quantum group $\G$ to be discrete or compact.
\end{remarks}
\renewcommand{\baselinestretch}{1.0}
\dated
\renewcommand{\baselinestretch}{1.2}
\vfill
\begin{tabbing}
\textit{Address}: \= Department of Mathematical and Statistical Sciences \\
\> University of Alberta \\
\> Edmonton, Alberta \\
\> Canada T6G 2G1 \\[\medskipamount]
\textit{E-mail}: \> \texttt{vrunde@ualberta.ca} \\[\medskipamount]
\textit{URL}: \> \texttt{http://www.math.ualberta.ca/$^\sim$runde/}   
\end{tabbing} 
\end{document}